\documentclass[secthm,seceqn,amsthm,ussrhead,12pt]{amsart}
\usepackage[utf8]{inputenc}
\usepackage[english]{babel}
\usepackage{amssymb,amsmath,amsthm,amsfonts,xcolor,enumerate,hyperref,comment,longtable,cleveref}

\usepackage{times}
\usepackage{cite}
\usepackage{pdflscape}
\usepackage{ulem}
\usepackage[mathcal]{euscript}
\usepackage{tikz}
\usepackage{hyperref}
\usepackage{cancel}
\usepackage{stmaryrd}
\usepackage{leftidx}

\usetikzlibrary{arrows}

\newtheorem{theorem}{Theorem}
\newtheorem{lemma}[theorem]{Lemma}

\newtheorem{remark}[theorem]{Remark}

\newtheorem{definition}[theorem]{Definition}

\newenvironment{Proof}[1][Proof.]{\begin{trivlist}
\item[\hskip \labelsep {\bfseries #1}]}{\flushright
$\Box$\end{trivlist}}

\usepackage{stmaryrd}
\usepackage{xcolor}

\begin{document}
\noindent{\bf
The algebraic and geometric classification of nilpotent bicommutative algebras}
\footnote{
This work was supported by  
RFBR 18-31-00001; 
FAPESP  18/15712-0;
MTM2016-79661-P;
FPU scholarship (Spain).
}$^,$\footnote{Corresponding Author: kaygorodov.ivan@gmail.com  }

   \

   {\bf  
   Ivan   Kaygorodov$^{a}$,
   María Pilar Paez Guillán$^{b}$ \&
   Vasily Voronin$^{c}$}

\

 \

{\tiny

$^{a}$ CMCC, Universidade Federal do ABC, Santo Andr\'e, Brazil.

$^{b}$ University of Santiago de Compostela, Santiago de Compostela, Spain.

$^{c}$ Novosibirsk State University, Novosibirsk, Russia.

\smallskip

   E-mail addresses:

\smallskip
Ivan   Kaygorodov (kaygorodov.ivan@gmail.com)

María Pilar Paez Guillán (pilar.paez@usc.es)

Vasily Voronin (voronin.vasily@gmail.com)

}

\

\

\noindent{\bf Abstract}:
{\it We classify the $4$-dimensional nilpotent bicommutative
 algebras over $\mathbb C$ from both algebraic and geometric approaches.}

\


\


\section*{Introduction}

One of the classical problems in the theory of non-associative algebras is to classify (up to isomorphism) the algebras of dimension $n$ from a certain variety defined by some family of polynomial identities. It is typical to focus on small dimensions, and there are two main directions for the classification: algebraic and geometric. Varieties as Jordan, Lie, Leibniz or Zinbiel algebras have been studied from these two approaches (\!\cite{ack, cfk19, gkks, degr3, usefi1, degr1, degr2, ha16, hac18, kv16} and \cite{casas, maria,bb14, BC99, cfk19, GRH, GRH2,gkks, ikv17, ikv18, kppv, kpv, kv16, S90}, respectively).
In the present paper, we give the algebraic  and geometric classification of
$4$-dimensional nilpotent bicommutative  algebras.

The variety of bicommutative algebras is defined by the following identities of right- and left-commutativity:
\[
\begin{array}{rclllrcl}
(xy)z &=& (xz)y, & \ &  x(yz) &=& y(xz).
\end{array} \]
It admits the commutative associative algebras as a subvariety. One-sided commutative algebras first appeared in the paper by Cayley \cite{cayley} in 1857. The structure of the free bicommutative algebra of countable rank and its main numerical invariants were described by
Dzhumadildaev, Ismailov, and Tulenbaev \cite{dit11}, see also the announcement \cite{dt03}. Bicommutative algebras were also studied in \cite{drensky1, drensky2, di}.

The key step in our method for algebraically classifying bicommutative nilpotent algebras is the calculation of central extensions of smaller algebras. \!In the theory of Lie groups, Lie algebras and their representations, the Lie algebra extensions
are enlargements of a given Lie algebra $\mathfrak{g}$ by another Lie algebra $\mathfrak{h}$, and they may
arise in several ways: for instance, as a direct sum of
two Lie algebras (trivial extensions), or when constructing a Lie algebra from projective group representations. Besides the simple trivial extensions, there are more important types, such as split and central extensions. 
There are important algebras which are central extensions: the Virasoro
algebra is the universal central extension of the Witt algebra, and the Heisenberg algebra is
the central extension of a commutative Lie algebra  \cite[Chapter 18]{bkk}, for example.
Further indicators of the importance of central extensions are the following facts. A central extension together with an extension by a derivation of a polynomial loop algebra
over a finite-dimensional simple Lie algebra give a Lie algebra isomorphic to a
non-twisted affine Kac--Moody algebra \cite[Chapter 19]{bkk}. Also, using the centrally extended loop
algebra, we can construct a current algebra in two spacetime dimensions.

Yet, central extensions not only play an important role in mathematics but also in physics, especially in different areas of quantum theory. For example, in quantum mechanics they first appeared in Wigner\'{}s theorem, which states that a symmetry of a quantum
mechanical system determines an (anti-)unitary transformation of a Hilbert space.
Also, we find them in the quantum theory
of conserved currents of a Lagrangian. These currents span an algebra which is closely
related to the universal central extensions
of loop algebras, namely affine Kac--Moody algebras.
Roughly speaking, central extensions are needed in physics because the symmetry group of a quantized
system is usually a central extension of the classical symmetry group, and in the same way
the corresponding symmetry Lie algebra of the quantum system is, in general, a central
extension of the classical symmetry algebra. In particular, Kac--Moody algebras have been conjectured
to be the symmetry groups of a unified superstring theory. Also, the centrally extended Lie
algebras play a dominant role in quantum field theory, particularly in conformal field
theory, string theory and in $M$-theory.

With this background, it comes as no surprise that the central extensions of Lie and non-Lie algebras have been exhaustively studied for years. It is interesting both to describe them and to use them to classify different varieties of algebras \cite{omirov,ha17,hac16,kkl18,ss78,zusmanovich}.
Firstly, Skjelbred and Sund devised a method for classifying nilpotent Lie algebras employing central extensions  \cite{ss78}.
Using this method, all the non-Lie central extensions of  all $4$-dimensional Malcev algebras were described  afterwards \cite{hac16}, and also
all the non-associative central extensions of $3$-dimensional Jordan algebras \cite{ha17},
all the anticommutative central extensions of the $3$-dimensional anticommutative algebras \cite{cfk182},
and all the central extensions of the $2$-dimensional algebras \cite{cfk18}.
Moreover, the method is especially indicated for the classification of nilpotent algebras
(see, for example, \cite{ha16n}),
and it was used to describe
all the $4$-dimensional nilpotent associative algebras \cite{degr1},
all the $4$-dimensional nilpotent Novikov algebras \cite{kkk18},
all the $5$-dimensional nilpotent Jordan algebras \cite{ha16},
all the $5$-dimensional nilpotent restricted Lie algebras \cite{usefi1},
all the $6$-dimensional nilpotent Lie algebras \cite{degr3,degr2},
all the $6$-dimensional nilpotent Malcev algebras \cite{hac18}
and some others.

\section{The algebraic classification of nilpotent bicommutative algebras}

\subsection{Method of classification of nilpotent algebras}

The objective of this section is to give an analogue of the Skjelbred-Sund method for classifying nilpotent bicommutative algebras. As other analogues of this method were carefully explained in, for example, \cite{ha17,hac16,cfk18}, we will give only some important definitions, and refer the interested reader to the previous sources. We will also employ their notations.

Let $({\bf A}, \cdot)$ be a bicommutative  algebra over  $\mathbb C$
and $\mathbb V$ a vector space over ${\mathbb C}$. We define the $\mathbb C$-linear space $Z^{2}\left(
\bf A,\mathbb V \right) $  as the set of all  bilinear maps $\theta  \colon {\bf A} \times {\bf A} \longrightarrow {\mathbb V}$
such that

\[ \theta(xy,z)=\theta(xz,y), \]
\[ \theta(x,yz)= \theta(y,xz). \]
These maps will be called \emph{cocycles}. Consider a
linear map $f$ from $\bf A$ to  $\mathbb V$, and set $\delta f\colon {\bf A} \times
{\bf A} \longrightarrow {\mathbb V}$ with $\delta f  (x,y ) =f(xy )$. Then, $\delta f$ is a cocycle, and we define $B^{2}\left(
{\bf A},{\mathbb V}\right) =\left\{ \theta =\delta f\ : f\in Hom\left( {\bf A},{\mathbb V}\right) \right\} $, a linear subspace of $Z^{2}\left( {\bf A},{\mathbb V}\right) $; its elements are called
\emph{coboundaries}. The \emph{second cohomology space} $H^{2}\left( {\bf A},{\mathbb V}\right) $ is defined to be the quotient space $Z^{2}
\left( {\bf A},{\mathbb V}\right) \big/B^{2}\left( {\bf A},{\mathbb V}\right) $.

\

Let $\operatorname{Aut}({\bf A}) $ be the automorphism group of the bicommutative algebra ${\bf A} $ and let $\phi \in \operatorname{Aut}({\bf A})$. Every $\theta \in
Z^{2}\left( {\bf A},{\mathbb V}\right) $ defines $\phi \theta (x,y)
=\theta \left( \phi \left( x\right) ,\phi \left( y\right) \right) $, with $\phi \theta \in Z^{2}\left( {\bf A},{\mathbb V}\right) $. It is easily checked that $\operatorname{Aut}({\bf A})$
acts on $Z^{2}\left( {\bf A},{\mathbb V}\right) $, and that
 $B^{2}\left( {\bf A},{\mathbb V}\right) $ is invariant under the action of $\operatorname{Aut}({\bf A}).$  
 So, we have that $\operatorname{Aut}({\bf A})$ acts on $H^{2}\left( {\bf A},{\mathbb V}\right)$.

\

Let $\bf A$ be a bicommutative  algebra of dimension $m<n$ over  $\mathbb C$, ${\mathbb V}$ a $\mathbb C$-vector
space of dimension $n-m$ and $\theta$ a cocycle, and consider the direct sum ${\bf A}_{\theta } = {\bf A}\oplus {\mathbb V}$ with the
bilinear product `` $\left[ -,-\right] _{{\bf A}_{\theta }}$'' defined by $\left[ x+x^{\prime },y+y^{\prime }\right] _{{\bf A}_{\theta }}=
 xy +\theta(x,y) $ for all $x,y\in {\bf A},x^{\prime },y^{\prime }\in {\mathbb V}$.
It is straightforward that ${\bf A_{\theta}}$ is a bicommutative algebra if and only if $\theta \in Z^2({\bf A}, {\mathbb V})$; it is  called an $(n-m)$-\emph{dimensional central extension} of ${\bf A}$ by ${\mathbb V}$.

We also call the
set $\operatorname{Ann}(\theta)=\left\{ x\in {\bf A}:\theta \left( x, {\bf A} \right)+ \theta \left({\bf A} ,x\right) =0\right\} $
the \emph{annihilator} of $\theta $. We recall that the \emph{annihilator} of an  algebra ${\bf A}$ is defined as
the ideal $\operatorname{Ann}(  {\bf A} ) =\left\{ x\in {\bf A}:  x{\bf A}+ {\bf A}x =0\right\}$. Observe
 that
$\operatorname{Ann}\left( {\bf A}_{\theta }\right) =\big(\operatorname{Ann}(\theta) \cap\operatorname{Ann}({\bf A})\big)
 \oplus {\mathbb V}$.

\

\begin{definition}
Given an algebra ${\bf A}$, if ${\bf A}=I\oplus \mathbb Cx$
is a direct sum   for some $x\in\operatorname{Ann}({\bf A})$, then $\mathbb Cx$ is called an \emph{annihilator component} of ${\bf A}$.
\end{definition}
\begin{definition}
A central extension of an algebra $\bf A$ without annihilator component is called a \emph{non-split central extension}.
\end{definition}

\

The following result is fundamental for the classification method.

\begin{lemma}
Let ${\bf A}$ be an $n$-dimensional bicommutative algebra such that $\dim(\operatorname{Ann}({\bf A}))=m\neq0$. Then there exists, up to isomorphism, a unique $(n-m)$-dimensional bicommutative  algebra ${\bf A}'$ and a bilinear map $\theta \in Z^2({\bf A}, {\mathbb V})$ with $\operatorname{Ann}({\bf A})\cap\operatorname{Ann}(\theta)=0$, where $\mathbb V$ is a vector space of dimension m, such that ${\bf A} \cong {{\bf A}'}_{\theta}$ and
 ${\bf A}/\operatorname{Ann}({\bf A})\cong {\bf A}'$.
\end{lemma}

For the proof, we refer the reader to~\cite[Lemma 5]{hac16}.


\

Now, we seek a condition on the cocycles to know when two $(n-m)$-central extensions are isomorphic.
Let us fix a basis $e_{1},\ldots ,e_{s}$ of ${\mathbb V}$, and $
\theta \in Z^{2}\left( {\bf A},{\mathbb V}\right) $. Then $\theta $ can be uniquely
written as $\theta \left( x,y\right) =
\displaystyle \sum_{i=1}^{s} \theta _{i}\left( x,y\right) e_{i}$, where $\theta _{i}\in
Z^{2}\left( {\bf A},\mathbb C\right) $. It holds that $\theta \in
B^{2}\left( {\bf A},{\mathbb V}\right) $\ if and only if all $\theta _{i}\in B^{2}\left( {\bf A},
\mathbb C\right) $, and it also holds that $\operatorname{Ann}(\theta)=\operatorname{Ann}(\theta _{1})\cap\operatorname{Ann}(\theta _{2})\ldots \cap\operatorname{Ann}(\theta _{s})$. 
Furthermore, if $\operatorname{Ann}(\theta)\cap \operatorname{Ann}\left( {\bf A}\right) =0$, then ${\bf A}_{\theta }$ has an
annihilator component if and only if $\left[ \theta _{1}\right] ,\left[
\theta _{2}\right] ,\ldots ,\left[ \theta _{s}\right] $ are linearly
dependent in $H^{2}\left( {\bf A},\mathbb C\right)$ (see \cite[Lemma 13]{hac16}).

\;

Recall that, given a finite-dimensional vector space ${\mathbb V}$ over $\mathbb C$, the \emph{Grassmannian} $G_{k}\left( {\mathbb V}\right) $ is the set of all $k$-dimensional
linear subspaces of $ {\mathbb V}$. Let $G_{s}\left( H^{2}\left( {\bf A},\mathbb C\right) \right) $ be the Grassmannian of subspaces of dimension $s$ in
$H^{2}\left( {\bf A},\mathbb C\right) $.
 For $W=\left\langle
\left[ \theta _{1}\right] ,\left[ \theta _{2}\right] ,\dots,\left[ \theta _{s}
\right] \right\rangle \in G_{s}\left( H^{2}\left( {\bf A},\mathbb C
\right) \right) $ and $\phi \in \operatorname{Aut}({\bf A})$, define $\phi W=\left\langle \left[ \phi \theta _{1}\right]
,\left[ \phi \theta _{2}\right] ,\dots,\left[ \phi \theta _{s}\right]
\right\rangle $. It holds that $\phi W\in G_{s}\left( H^{2}\left( {\bf A},\mathbb C \right) \right) $, and this induces an action of $\operatorname{Aut}({\bf A})$ on $G_{s}\left( H^{2}\left( {\bf A},\mathbb C\right) \right) $. We denote the orbit of $W\in G_{s}\left(
H^{2}\left( {\bf A},\mathbb C\right) \right) $ under this action  by $\operatorname{Orb}(W)$. Let
\[
W_{1}=\left\langle \left[ \theta _{1}\right] ,\left[ \theta _{2}\right] ,\dots,
\left[ \theta _{s}\right] \right\rangle ,W_{2}=\left\langle \left[ \vartheta
_{1}\right] ,\left[ \vartheta _{2}\right] ,\dots,\left[ \vartheta _{s}\right]
\right\rangle \in G_{s}\left( H^{2}\left( {\bf A},\mathbb C\right)
\right).
\]
Similarly to~\cite[Lemma 15]{hac16}, in case $W_{1}=W_{2}$, it holds that \[ \bigcap\limits_{i=1}^{s}\operatorname{Ann}(\theta _{i})\cap \operatorname{Ann}\left( {\bf A}\right) = \bigcap\limits_{i=1}^{s}
\operatorname{Ann}(\vartheta _{i})\cap\operatorname{Ann}( {\bf A}) ,\] 
and therefore the set
\[
T_{s}({\bf A}) =\left\{ W=\left\langle \left[ \theta _{1}\right] ,
\left[ \theta _{2}\right] ,\dots,\left[ \theta _{s}\right] \right\rangle \in
G_{s}\left( H^{2}\left( {\bf A},\mathbb C\right) \right) : \bigcap\limits_{i=1}^{s}\operatorname{Ann}(\theta _{i})\cap\operatorname{Ann}({\bf A}) =0\right\}
\]
is well defined, and it is also stable under the action of $\operatorname{Aut}({\bf A})$ (see~\cite[Lemma 16]{hac16}).

\

Now, let ${\mathbb V}$ be an $s$-dimensional linear space and let us denote by
$E\left( {\bf A},{\mathbb V}\right) $ the set of all non-split $s$-dimensional central extensions of ${\bf A}$ by
${\mathbb V}$. We can write
\[
E\left( {\bf A},{\mathbb V}\right) =\left\{ {\bf A}_{\theta }:\theta \left( x,y\right) = \sum_{i=1}^{s}\theta _{i}\left( x,y\right) e_{i} \ \ \text{and} \ \ \left\langle \left[ \theta _{1}\right] ,\left[ \theta _{2}\right] ,\dots,
\left[ \theta _{s}\right] \right\rangle \in T_{s}({\bf A}) \right\} .
\]

Finally, we are prepared to state our main result, which can be proved as \cite[Lemma 17]{hac16}.

\begin{lemma}
 Let ${\bf A}_{\theta },{\bf A}_{\vartheta }\in E\left( {\bf A},{\mathbb V}\right) $. Suppose that $\theta \left( x,y\right) =  \displaystyle \sum_{i=1}^{s}
\theta _{i}\left( x,y\right) e_{i}$ and $\vartheta \left( x,y\right) =
\displaystyle \sum_{i=1}^{s} \vartheta _{i}\left( x,y\right) e_{i}$.
Then the bicommutative algebras ${\bf A}_{\theta }$ and ${\bf A}_{\vartheta } $ are isomorphic
if and only if
$$\operatorname{Orb}\left\langle \left[ \theta _{1}\right] ,
\left[ \theta _{2}\right] ,\dots,\left[ \theta _{s}\right] \right\rangle =
\operatorname{Orb}\left\langle \left[ \vartheta _{1}\right] ,\left[ \vartheta
_{2}\right] ,\dots,\left[ \vartheta _{s}\right] \right\rangle .$$
\end{lemma}

Then, it exists a bijective correspondence between the set of $\operatorname{Aut}({\bf A})$-orbits on $T_{s}\left( {\bf A}\right) $ and the set of
isomorphism classes of $E\left( {\bf A},{\mathbb V}\right) $. Consequently we have a
procedure that allows us, given a bicommutative algebra ${\bf A}'$ of
dimension $n-s$, to construct all non-split central extensions of ${\bf A}'$.

\; \;

{\centerline {\textsl{Procedure}}}

Let ${\bf A}'$ be a bicommutative algebra of dimension $n-s $.

\begin{enumerate}
\item Determine $H^{2}( {\bf A}',\mathbb {C}) $, $\operatorname{Ann}({\bf A}')$ and $\operatorname{Aut}({\bf A}')$.

\item Determine the set of $\operatorname{Aut}({\bf A}')$-orbits on $T_{s}({\bf A}') $.

\item For each orbit, construct the bicommutative algebra associated with a
representative of it.
\end{enumerate}

\

\subsection{Notations}
Let ${\bf A}$ be a bicommutative algebra and fix
a basis $e_{1},e_{2},\dots,e_{n}$. We define the bilinear form
$\Delta _{ij} \colon {\bf A}\times {\bf A}\longrightarrow \mathbb C$
by $\Delta _{ij}\left( e_{l},e_{m}\right) = \delta_{il}\delta_{jm}$.
Then the set $\left\{ \Delta_{ij}:1\leq i, j\leq n\right\} $ is a basis for the linear space of
the bilinear forms on ${\bf A}$, and in particular, every $\theta \in
Z^{2}\left( {\bf A},\mathbb V\right) $ can be uniquely written as $
\theta = \displaystyle \sum_{1\leq i,j\leq n} c_{ij}\Delta _{{i}{j}}$, where $
c_{ij}\in \mathbb C$.
Let us fix the following notations:

$$\begin{array}{lll}
{\mathcal B}^{i*}_j& \mbox{---}& j\mbox{th }i\mbox{-dimensional nilpotent ``non-pure'' bicommutative algebra (with identity $xyz=0$}); \\
{\mathcal B}^i_j& \mbox{---}& j\mbox{th }i\mbox{-dimensional nilpotent ``pure'' bicommutative algebra (without identity $xyz=0$)}; \\
{\mathfrak{N}}_i& \mbox{---}& i\mbox{-dimensional algebra with zero product}; \\
({\bf A})_{i,j}& \mbox{---}& j\mbox{th }i\mbox{-dimensional central extension of }\bf A. \\
\end{array}$$

\subsection{The algebraic classification of  $3$-dimensional nilpotent bicommutative algebras}
There are no nontrivial $1$-dimensional nilpotent bicommutative algebras, and 
there is only one nontrivial $2$-dimensional nilpotent bicommutative algebra
(namely, the non-split central extension of the $1$-dimensional algebra with zero product):

$$\begin{array}{ll llll}
{\mathcal B}^{2*}_{01} &:& (\mathfrak{N}_1)_{2,1} &:& e_1 e_1 = e_2.\\
\end{array}$$

From this algebra, we construct the $3$-dimensional nilpotent bicommutative algebra ${\mathcal B}^{3*}_{01}={\mathcal B}^{2*}_{01}\oplus{\mathbb C e_3}.$

Also, the reference \cite{cfk18} gives the description of all central extensions of  ${\mathcal B}^{2*}_{01}$ and $\mathfrak{N}_2$.
Choosing the bicommutative algebras between them,
we have the classification of all non-split $3$-dimensional nilpotent bicommutative algebras:

$$\begin{array}{ll llllllllllll}
{\mathcal B}^{3*}_{02} &:& (\mathfrak{N}_2)_{3,1} &:& e_1 e_1 = e_3, &  e_2 e_2=e_3; \\
{\mathcal B}^{3*}_{03} &:& (\mathfrak{N}_2)_{3,2} &:& e_1 e_2=e_3, & e_2 e_1=-e_3;   \\
{\mathcal B}^{3*}_{04}(\alpha)_{\alpha \neq 0} &:& (\mathfrak{N}_2)_{3,3} &:&
e_1 e_1 = \alpha e_3,  & e_2 e_1=e_3,  & e_2 e_2=e_3; &  \\
{\mathcal B}^{3*}_{04}(0) &:& (\mathfrak{N}_2)_{3,3} &:&   e_1 e_2=e_3;  \\
{\mathcal B}^3_{01} &:& ({\mathcal B}^{2*}_{01} )_{3,1} &:& e_1 e_1 = e_2,  & e_2 e_1=e_3;  \\
{\mathcal B}^3_{02}(\alpha) &:& ({\mathcal B}^{2*}_{01} )_{3,2} &:& e_1 e_1 = e_2,  & e_1 e_2=e_3,  & e_2 e_1 = \alpha e_3.  \\
\end{array} $$

\subsection{$1$-dimensional central extensions of $3$-dimensional  nilpotent bicommutative algebras}
\label{centrext}
\subsubsection{The description of second cohomology space of  $3$-dimensional nilpotent bicommutative algebras}

\
In the following table we give the description of the second cohomology space of  $3$-dimensional nilpotent bicommutative algebras.

{\tiny
$$
\begin{array}{|l|l|l|l|}
\hline
\bf A  & Z^{2}\left( {\bf A}\right)  & B^2({\bf A}) & H^2({\bf A}) \\
\hline
\hline
{\mathcal B}^{3*}_{01} &  \langle
\Delta_{11},\Delta_{12}, \Delta_{13}, \Delta_{21},\Delta_{31}, \Delta_{33}\rangle
&\langle \Delta_{11} \rangle&
\langle[\Delta_{12}], [\Delta_{13}], [\Delta_{21}], [\Delta_{31}], [\Delta_{33}] \rangle\\
\hline

{\mathcal B}^{3*}_{02} &  \langle  \Delta_{11},\Delta_{12}, \Delta_{21}, \Delta_{22} \rangle
& \langle\Delta_{11}+\Delta_{22}\rangle &  \langle [\Delta_{12}], [\Delta_{21}], [\Delta_{22}] \rangle \\
\hline
{\mathcal B}^{3*}_{03} & \langle  \Delta_{11},\Delta_{12}, \Delta_{21}, \Delta_{22} \rangle
& \langle\Delta_{12} -\Delta_{21} \rangle&  \langle [\Delta_{11}], [\Delta_{21}], [\Delta_{22}] \rangle \\
\hline
{\mathcal B}^{3*}_{04}(\alpha)_{\alpha\neq 0} & \langle  \Delta_{11},\Delta_{12}, \Delta_{21}, \Delta_{22} \rangle & \langle\alpha\Delta_{11}+\Delta_{21}+\Delta_{22}\rangle &  \langle [\Delta_{12}], [\Delta_{21}], [\Delta_{22}] \rangle \\

\hline
{\mathcal B}^{3*}_{04}(0) &
\langle \Delta_{11}, \Delta_{12}, \Delta_{13}, \Delta_{21}, \Delta_{22}, \Delta_{32} \rangle &
\langle \Delta_{12} \rangle &
\langle [\Delta_{11}], [\Delta_{13}], [\Delta_{21}], [\Delta_{22}], [\Delta_{32}] \rangle

\\

\hline
{\mathcal B}^3_{01} &\langle \Delta_{11},\Delta_{12}, \Delta_{21},\Delta_{31}\rangle &
\langle \Delta_{11}, \Delta_{21}\rangle &
\langle [\Delta_{12}], [\Delta_{31}] \rangle \\

\hline
{\mathcal B}^3_{02}(\alpha) &
\langle \Delta_{11},\Delta_{12}, \Delta_{21}, \alpha\Delta_{22}+\Delta_{13}+\alpha\Delta_{31} \rangle &
\langle \Delta_{11}, \Delta_{12} + \alpha \Delta_{21} \rangle &
\langle [\Delta_{21}], \alpha[\Delta_{22}]+[\Delta_{13}]+\alpha[\Delta_{31}] \rangle \\


\hline
\end{array}$$
}

\begin{remark}
From the description of the cocycles of the algebras ${\mathcal B}^{3*}_{02}$, ${\mathcal B}^{3*}_{03}$ and ${\mathcal B}^{3*}_{04}(\alpha)_{\alpha\neq 0}$, it follows that the 1-dimensional central extensions of these algebras are 2-dimensional central extensions of 2-dimensional nilpotent bicommutative algebras. Thanks to \cite{cfk18} we have the description of all non-split 2-dimensional central extensions of 2-dimensional nilpotent bicommutative algebras:

$$\begin{array}{ll lllll}
{\mathcal B}^{4}_{03} &:& ({\mathcal B}^{2*}_{01})_{4,1} &:& e_1 e_1 = e_2, & e_1 e_2 = e_4, & e_2 e_1 = e_3.\\
\end{array}$$
Then, in the following subsections we study the central extensions of the other algebras.
\end{remark}

\subsubsection{Central extensions of ${\mathcal B}^{3*}_{01}$}

Since the second cohomology spaces and automorphism groups of ${\mathcal B}^{3*}_{01}$ and ${\mathcal N}^{3*}_{01}$  (from \cite{kkk18}) coincide,  these algebras have the same central extensions.
Therefore, thanks to \cite{kkk18} we have all the new $4$-dimensional nilpotent bicommutative algebras constructed from   ${\mathcal B}^{3*}_{01}$:
\[ {\mathcal B}^4_{04}(\alpha), \ {\mathcal B}^4_{05}, \ {\mathcal B}^4_{06}(\alpha)_{\alpha\neq 0}, \ {\mathcal B}^4_{07}, \ {\mathcal B}^4_{08}, \ {\mathcal B}^4_{09}. \]
The multiplication tables of these algebras can be found in Appendix A.

\subsubsection{Central extensions of ${\mathcal B}^{3*}_{04}(0)$}

Let us use the following notations:
\[\nabla_1=[\Delta_{11}], \nabla_2=[\Delta_{13}], \nabla_3=[\Delta_{21}], \nabla_4=[\Delta_{22}],
\nabla_5=[\Delta_{32}].\]

The automorphism group of ${\mathcal B}^{3*}_{04}(0)$ consists of invertible matrices of the form

\[\phi=\left(
                             \begin{array}{ccc}
                               x & 0 & 0   \\
                               0 & y & 0  \\
                               z & t & xy                               \end{array}\right)
                               .\]

Since
\[
\phi^T
                           \left(\begin{array}{ccc}
                                \alpha_1& 0 & \alpha_2  \\
                                 \alpha_3 & \alpha_4 & 0  \\
                                 0 & \alpha_5 & 0 \\
                             \end{array}
                           \right)\phi
                           =\left(\begin{array}{ccc}
                                 x(x\alpha_1+z\alpha_2) & \alpha^* & x^2y\alpha_2   \\
                                 xy\alpha_3 & y(y\alpha_4+t\alpha_5)& 0  \\
                                 0 & xy^2\alpha_5 & 0 \\
                             \end{array}\right),\]
we have that the action of $\operatorname{Aut} ({\mathcal B}^{3*}_{04}(0))$ on the subspace
$\langle  \sum\limits_{i=1}^5\alpha_i \nabla_i \rangle$
is given by
$\langle  \sum\limits_{i=1}^5\alpha^*_i \nabla_i \rangle,$ where

\[
\begin{array}{rcl}
\alpha^*_1&=&x(x\alpha_1+z\alpha_2);\\
\alpha^*_2&=&x^2y\alpha_2;\\
\alpha^*_3&=&xy\alpha_3;\\
\alpha^*_4&=&y(y\alpha_4+t\alpha_5);\\
\alpha^*_5&=&xy^2\alpha_5.
\end{array}\]

It is easy to see that the elements $\alpha_1 \nabla_1 +\alpha_3\nabla_3 +\alpha_4\nabla_4$ give algebras
which are central extensions of $2$-dimensional algebras.
We find the following new cases:

\begin{enumerate}

\item $\alpha_2\neq0, \alpha_3\neq 0, \alpha_5\neq0,$ then choosing $x=\frac{\alpha_3}{\alpha_2}$,
      $y=\frac{\alpha_3}{\alpha_5}$, $z=-\frac{x\alpha_1}{\alpha_2}$ and
      $t=-\frac{y\alpha_4}{\alpha_5}$,
      we have the representative $\langle \nabla_2+\nabla_3+\nabla_5 \rangle$.

\item $\alpha_2 \neq 0, \alpha_3 = 0, \alpha_5 \neq 0,$ then choosing $y=\frac{x\alpha_2}{\alpha_5}$,
      $z=-\frac{x\alpha_1}{\alpha_2}$ and $t=-\frac{y\alpha_4}{\alpha_5}$,
      we have the representative $\langle \nabla_2+\nabla_5 \rangle$.

\item $\alpha_2 = 0, \alpha_3 \neq 0, \alpha_5 \neq 0,$ then:
\begin{enumerate}
    \item if $\alpha_1 \neq 0$, then choosing $y=\frac{\alpha_3}{\alpha_5}$,
    $x=\frac{y\alpha_3}{\alpha_1}$ and $t=-\frac{y\alpha_4}{\alpha_5}$,
    we have the representative $\langle \nabla_1+\nabla_3+\nabla_5 \rangle$.

    \item if $\alpha_1 = 0$, then choosing $y=\frac{\alpha_3}{\alpha_5}$ and
    $t=-\frac{y\alpha_4}{\alpha_5}$,
    we have the representative $\langle \nabla_3+\nabla_5 \rangle$.
\end{enumerate}

\item $\alpha_2 \neq 0, \alpha_3 \neq 0, \alpha_5 = 0,$ then:
\begin{enumerate}
    \item if $\alpha_4 \neq 0$, then choosing $x=\frac{\alpha_3}{\alpha_2}$,
    $y=\frac{x\alpha_3}{\alpha_4}$ and $z=-\frac{x\alpha_1}{\alpha_2}$,
    we have the representative $\langle \nabla_2+\nabla_3+\nabla_4 \rangle$.

    \item if $\alpha_4 = 0$, then choosing $x=\frac{\alpha_3}{\alpha_2}$
    and  $z=-\frac{x\alpha_1}{\alpha_2}$,
    we have the representative $\langle \nabla_2+\nabla_3 \rangle$.
\end{enumerate}

\item $\alpha_2 \neq 0, \alpha_3 = 0, \alpha_5 = 0$, then:
\begin{enumerate}
    \item if $\alpha_4 \neq 0$, then choosing $y=\frac{x^2\alpha_2}{\alpha_4}$ and
    $z=-\frac{x\alpha_1}{\alpha_2}$,
    we have the representative $\langle \nabla_2+\nabla_4 \rangle$.

    \item if $\alpha_4 = 0$, then choosing $z=-\frac{x\alpha_1}{\alpha_2}$,
    we have the representative $\langle \nabla_2 \rangle$.
\end{enumerate}

\item $\alpha_2 = 0, \alpha_3 = 0, \alpha_5 \neq 0$, then:
\begin{enumerate}
    \item  if $\alpha_1 \neq 0$, then choosing $x=\frac{y^2\alpha_5}{\alpha_1}$ and
    $t=-\frac{y\alpha_4}{\alpha_5}$,
    we have the representative $\langle \nabla_1+\nabla_5 \rangle$.

    \item if  $\alpha_1 = 0$ then choosing $t=-\frac{y\alpha_4}{\alpha_5}$,
    we have the representative $\langle \nabla_5 \rangle$.
\end{enumerate}

\end{enumerate}

Now we have all the new $4$-dimensional nilpotent bicommutative algebras constructed from  ${\mathcal
B}^{3*}_{04}(0):$
\[ {\mathcal B}^4_{10}, \ \ldots, {\mathcal B}^{4}_{19}. \]
The multiplication tables of these algebras can be found in Appendix A.

\subsubsection{Central extensions of ${\mathcal B}^{3}_{01}$}

Let us use the following notations:
\[\nabla_1=[\Delta_{12}], \nabla_2=[\Delta_{31}].\]

The automorphism group of ${\mathcal B}^{3}_{01}$ consists of invertible matrices of the form

\[\phi=\left(
                             \begin{array}{ccc}
                               x & 0 & 0   \\
                               y & x^2 & 0  \\
                               z & xy & x^3                               \end{array}\right)
                               .\]

Since
\[
\phi^T
                           \left(\begin{array}{ccc}
                                0& \alpha_1 & 0  \\
                                 0 & 0 & 0  \\
                                 \alpha_2 & 0 & 0 \\
                             \end{array}
                           \right)\phi
                           =\left(\begin{array}{ccc}
                                 \alpha^* & x^3\alpha_1 & 0   \\
                                 \alpha^{**}& 0 & 0  \\
                                 x^4\alpha_2 & 0 & 0 \\
                             \end{array}\right),\]
we have that the action of $\operatorname{Aut} ({\mathcal B}^{3}_{01})$ on the subspace
$\langle  \sum\limits_{i=1}^2\alpha_i \nabla_i \rangle$
is given by
$\langle  \sum\limits_{i=1}^2\alpha^*_i \nabla_i \rangle,$ where

\[
\begin{array}{rcl}
\alpha^*_1&=&x^3\alpha_1;\\
\alpha^*_2&=&x^4\alpha_2.\\
\end{array}\]

It is straightforward that the elements $\alpha_1 \nabla_1$ lead to central extensions of $2$-dimensional algebras.
The new cases are following:

\begin{enumerate}

\item $\alpha_1\neq0, \alpha_2\neq 0.$ Choosing $x=\frac{\alpha_1}{\alpha_2}$,
      we have the representative $\langle \nabla_1+\nabla_2 \rangle$.

\item $\alpha_1 = 0, \alpha_2 \neq 0.$  Choosing $x=\frac{1}{\sqrt[4]{\alpha_2}}$,
      we have the representative $\langle \nabla_2\rangle$.
\end{enumerate}

Now we have all the new $4$-dimensional nilpotent bicommutative algebras constructed from  ${\mathcal
B}^{3}_{01}:$
\[ {\mathcal B}^4_{20},\ {\mathcal B}^{4}_{21}. \]
The multiplication tables of these algebras can be found in Appendix A.

\subsubsection{Central extensions of ${\mathcal B}^{3}_{02}(\alpha)$}

Let us use the following notations:
\[\nabla_1=[\Delta_{21}], \nabla_2=\alpha[\Delta_{22}]+[\Delta_{13}]+\alpha[\Delta_{31}].\]

The automorphism group of ${\mathcal B}^{3}_{02}(\alpha)$ consists of invertible matrices of the form

\[\phi=\left(
                             \begin{array}{ccc}
                               x & 0 & 0   \\
                               y & x^2 & 0  \\
                               z & (\alpha+1)xy & x^3                               \end{array}\right)
                               .\]

Since
\[
\phi^T
                           \left(\begin{array}{ccc}
                                0& 0 & \alpha_2  \\
                                 \alpha_1 & \alpha\alpha_2 & 0  \\
                                 \alpha\alpha_2 & 0 & 0 \\
                             \end{array}
                           \right)\phi
                           =\left(\begin{array}{ccc}
                                 \alpha^* & \alpha^{**} & \alpha^*_2   \\
                                 \alpha_1^*+\alpha\alpha^{**} & \alpha \alpha^*_2 & 0  \\
                                 \alpha \alpha^*_2 & 0 & 0 \\
                             \end{array}\right),\]
we have that the action of $\operatorname{Aut} ({\mathcal B}^{3}_{02}(\alpha))$ on the subspace
$\langle  \sum\limits_{i=1}^2\alpha_i \nabla_i \rangle$
is given by
$\langle  \sum\limits_{i=1}^2\alpha^*_i \nabla_i \rangle,$ where

\[
\begin{array}{rcl}
\alpha^*_1&=&x^2(x\alpha_1+ \alpha(1-\alpha)y\alpha_2);\\
\alpha^*_2&=&x^4\alpha_2.\\
\end{array}\]

The element $\alpha_1 \nabla_1$ gives a central extension of a $2$-dimensional algebra, 
then we will consider only cases with $\alpha_2\neq 0.$
We find the following new cases:

\begin{enumerate}

\item $\alpha=0$ or $\alpha=1,$ then:

    \begin{enumerate}
        \item if $\alpha_1 \neq 0$, then choosing $x=\frac{\alpha_1}{\alpha_2}$,
             we have the representative $\langle \nabla_1+\nabla_2 \rangle$.

      \item if $\alpha_1 = 0$, then choosing $x=\frac{1}{\sqrt[4]{\alpha_2}}$,
        we have the representative $\langle \nabla_2 \rangle$.
    \end{enumerate}

\item $\alpha \neq 0,1,$ then choosing $x=\frac{1}{\sqrt[4]{\alpha_2}}$ and $y=\frac{-x\alpha_1}{\alpha_2\alpha(1-\alpha)}$, we have the representative $\langle \nabla_2\rangle$.
\end{enumerate}

Now we have all the new $4$-dimensional nilpotent bicommutative algebras constructed from  ${\mathcal
B}^{3}_{02}(\alpha):$
\[ {\mathcal B}^4_{22},\ {\mathcal B}^{4}_{23},\ {\mathcal B}^{4}_{24}(\alpha). \]
The multiplication tables of these algebras can be found in Appendix A.

\subsection{The algebraic classification of $4$-dimensional nilpotent bicommutative algebras}

We distinguish two main classes of bicommutative algebras: the ``pure'' and the ``non-pure'' ones. By the non-pure ones, we mean those satisfying the identities $(xy)z=0$ and $x(yz)=0$; the pure ones are the rest. 

These ``trivial'' algebras can be considered in many varieties of algebras defined by polynomial identities of degree $3$ (associative, Leibniz, Zinbiel\dots), and they can be expressed as central extensions of suitable algebras with zero product. Those with dimension $4$ are already classified:  the list of the non-anticommutative ones can be found in \cite{demir}, and there is only one  nilpotent and anticommutative.


Regarding the pure $4$-dimensional nilpotent bicommutative algebras, we have the following theorem, whose proof is  based on the classification of $3$-dimensional nilpotent bicommutative algebras and the results of Section \ref{centrext}.


\begin{theorem}
Let $\mathcal B$ be a nonzero  $4$-dimensional nilpotent pure bicommutative algebra over $\mathbb C$.
Then, $\mathcal B$ is isomorphic to one of the algebras listed in Table A (see Appendix).
\end{theorem}

\section{The geometric classification of nilpotent bicommutative algebras}

\subsection{Definitions and notation}
Let $\mathbb V$ be an $n$-dimensional vector space. The set $Hom(\mathbb V \otimes \mathbb V,\mathbb V) \cong \mathbb V^* \otimes \mathbb V^* \otimes \mathbb V$
is a vector space of dimension $n^3$, and it has the structure of the affine variety $\mathbb{C}^{n^3}.$ Indeed, if we fix a basis $e_1,\dots,e_n$ of $\mathbb V$, then any $\mu\in Hom(\mathbb V \otimes \mathbb V,\mathbb V)$ is determined by $n^3$ structure constants $c_{ij}^k\in\mathbb{C}$ such that
$\mu(e_i\otimes e_j)=\sum\limits_{k=1}^nc_{ij}^ke_k$. A subset of $Hom(\mathbb V \otimes \mathbb V,\mathbb V)$ is called {\it Zariski-closed} if it can be defined by a set of polynomial equations in the variables $c_{ij}^k$ ($1\le i,j,k\le n$).

Let $T$ be such a set of polynomial identities.
It holds that every algebra structure on $\mathbb V$ satisfying polynomial identities from $T$ forms a Zariski-closed subset of the variety $Hom(\mathbb V \otimes \mathbb V,\mathbb V)$; it is denoted by $\mathbb{L}(T)$.
There exists a natural action of the general linear group $GL(\mathbb V)$ on $\mathbb{L}(T)$ defined by
$$ (g * \mu )(x\otimes y) = g\mu(g^{-1}x\otimes g^{-1}y)$$
for $x,y\in \mathbb V$, $\mu\in \mathbb{L}(T)\subset Hom(\mathbb V \otimes\mathbb V, \mathbb V)$ and $g\in GL(\mathbb V)$.
Then, $\mathbb{L}(T)$ can be decomposed into $GL(\mathbb V)$-orbits corresponding to the isomorphism classes of the algebras.
We will denote by $O(\mu)$ the orbit of $\mu\in\mathbb{L}(T)$ under the action of $GL(\mathbb V)$, and by $\overline{O(\mu)}$ the Zariski closure of $O(\mu)$.

Let $\mathcal A$ and $\mathcal B$ be two $n$-dimensional algebras satisfying the identities from $T$, and let $\mu,\lambda \in \mathbb{L}(T)$ represent $\mathcal A$ and $\mathcal B$, respectively.
We say that $\mathcal A$ degenerates to $\mathcal B$, and write $\mathcal A\to \mathcal B$, if $\lambda\in\overline{O(\mu)}$.
Note that, in particular, it holds that $\overline{O(\lambda)}\subset\overline{O(\mu)}$. Hence, the definition of a degeneration does not depend on the choice of $\mu$ and $\lambda$. If $\mathcal A\not\cong \mathcal B$, then the assertion $\mathcal A\to \mathcal B$ is called a {\it proper degeneration}. Also, we write $\mathcal A\not\to \mathcal B$ if $\lambda\not\in\overline{O(\mu)}$.

Now consider $\mathcal A(*):=\{\mathcal A(\alpha)\}_{\alpha\in I}$ a family of algebras parameterized by $\alpha$,
and let $\mathcal A(\alpha)$, for $\alpha\in I$, be represented by the structure $\mu(\alpha)\in\mathbb{L}(T)$. 
Then $\mathcal A(*)\to \mathcal B$ means $\lambda\in\overline{\{O(\mu(\alpha))\}_{\alpha\in I}}$, and $\mathcal A(*)\not\to \mathcal B$ means $\lambda\not\in\overline{\{O(\mu(\alpha))\}_{\alpha\in I}}$.

Moreover, we call $\mathcal A$  {\it rigid} in $\mathbb{L}(T)$ if $O(\mu)$ is an open subset of $\mathbb{L}(T)$.
 Recall that a subset of a variety is called irreducible if it cannot be represented as a union of two non-trivial closed subsets, and that
 a maximal irreducible closed subset of a variety is called an {\it irreducible component}.
It is well known that any affine variety can be represented as a finite union of its irreducible components in a unique way.
Then, we have the following characterization of rigidity:
$\mathcal A$ is rigid in $\mathbb{L}(T)$ if and only if $\overline{O(\mu)}$ is an irreducible component of $\mathbb{L}(T)$.

Henceforth, given the spaces $U$ and $W$, we will write simply $U>W$ instead of $dim\,U>dim\,W$.

\subsection{Method of the description of  degenerations of algebras}

In the present work we use the methods applied to Lie algebras in \cite{BC99,GRH,GRH2,S90}.
Let $\mathfrak{Der}(\mathcal A)$ denote the Lie algebra of derivations of $\mathcal A$.
Our first and useful consideration is that 
if $\mathcal A\to \mathcal B$ and $\mathcal A\not\cong \mathcal B$, then $\mathfrak{Der}(\mathcal A)<\mathfrak{Der}(\mathcal B)$ and $\mathcal{A}^2\geq \mathcal{B}^2$. Then, we will compute the dimensions of algebras of derivations and will check the assertion $\mathcal A\to \mathcal B$ only for such $\mathcal A$ and $\mathcal B$ that $\mathfrak{Der}(\mathcal A)<\mathfrak{Der}(\mathcal B)$. Among them, we will calculate the dimension of the squares of the algebras and check $\mathcal A\to \mathcal B$ only for such $\mathcal A$ and $\mathcal B$ that $\mathcal{A}^2\geq \mathcal{B}^2$.

Now, we explain our method for proving degenerations. 
Let $\mathcal A$, $\mathcal A(*)$ and $\mathcal B$ be as in Subsection 2.1. 
Fixed a basis $e_1,\dots, e_n$  of $\mathbb  V$, let $c_{ij}^k$ ($1\le i,j,k\le n$) be the structure constants of $\lambda$ in this basis. On the one hand, if there exist $a_i^j(t)\in\mathbb{C}$ ($1\le i,j\le n$, $t\in\mathbb{C}^*$) such that $E_i^t=\sum\limits_{j=1}^na_i^j(t)e_j$ ($1\le i\le n$) form a basis of $\mathbb V$ for any $t\in\mathbb{C}^*$, and the structure constants of $\mu$ in the basis $E_1^t,\dots, E_n^t$ are such polynomials $c_{ij}^k(t)\in\mathbb{C}[t]$ that $c_{ij}^k(0)=c_{ij}^k$, then $\mathcal A\to \mathcal B$. 
In this case  $E_1^t,\dots, E_n^t$ is called a {\it parametrized basis} for $\mathcal A\to \mathcal B$.
On the other hand, if we construct $a_i^j:\mathbb{C}^*\to \mathbb{C}$ ($1\le i,j\le n$) and $f: \mathbb{C}^* \to I$ such that $E_i^t=\sum\limits_{j=1}^na_i^j(t)e_j$ ($1\le i\le n$) form a basis of $\mathbb V$ for any  $t\in\mathbb{C}^*$, and the structure constants of $\mu_{f(t)}$ in the basis $E_1^t,\dots, E_n^t$ are such polynomials $c_{ij}^k(t)\in\mathbb{C}[t]$ that $c_{ij}^k(0)=c_{ij}^k$, then $\mathcal A(*)\to \mathcal B$. 
In this case  $E_1^t,\dots, E_n^t$ and $f(t)$ are called a parametrized basis and a {\it parametrized index} for $\mathcal A(*)\to \mathcal B$, respectively.

\subsection{The geometric classification of $4$-dimensional nilpotent bicommutative algebras}
The main result of the present section is the following theorem.

\begin{theorem}\label{geobl}
The variety of $4$-dimensional nilpotent bicommutative algebras has two irreducible components
defined by the rigid algebra  $ {\mathcal B}^{4}_{10}$
and the infinite family of algebras ${\mathcal B}^4_{24}(\alpha)$.
\end{theorem}

\begin{Proof}

From the considerations about the dimension of derivations in the previous subsection, it follows that there are not nilpotent bicommutative algebras degenerating to $\mathcal{B}_{10}^4$. Also, the dimension of the square of $\mathcal{B}_{10}^4$ is $2$, and $\mathcal{B}_{24}^4(\alpha)$ has $3$-dimensional square, so it cannot degenerate from $\mathcal{B}_{10}^4$. Therefore, if we prove that these two algebras degenerate to the rest of the nilpotent bicommutative algebras, the theorem is proved.

Recall that the full description of the degeneration system of $4$-dimensional trivial bicommutative algebras was given in \cite{kppv}.
Using the cited result, we have that the variety of $4$-dimensional trivial bicommutative algebras has two irreducible components given by the two following
families of algebras:

$$\begin{array}{lllllll}
\mathfrak{N}_2(\alpha)  & e_1e_1 = e_3, &e_1e_2 = e_4,  &e_2e_1 = -\alpha e_3, &e_2e_2 = -e_4; \\

\mathfrak{N}_3(\alpha)  & e_1e_1 = e_4, &e_1e_2 = \alpha e_4,  &e_2e_1 = -\alpha e_4, &e_2e_2 = e_4,  &e_3e_3 = e_4.
\end{array}$$

The algebra $\mathcal{B}^{4}_{10}$ degenerates to both $\mathfrak{N}_2(\alpha)$ and $\mathfrak{N}_3(\alpha)$. We will explain in detail the degeneration ${\mathcal B}^{4}_{10} \to {\mathfrak N}_{3}(\alpha)_{\alpha \neq 0, \pm i}$; as for ${\mathcal B}^{4}_{10} \to {\mathfrak N}_{2}(\alpha)$, it is similar, but easier. It can be found in Table B (Appendix A).

Let us consider the parametric basis of ${\mathcal B}^4_{10}:$
$F^t_i= \sum\limits_{j=i}^4 \alpha_{ij}(t)e_j.$
The multiplication table in the new basis is given below:
$$\begin{array}{rclrcl}

F^t_1F^t_1 &=& \multicolumn{4}{l}{\frac{\alpha_{11}\alpha_{12}}{\alpha_{33}}F^t_3+\frac{\alpha_{11}\alpha_{13}+\alpha_{11}\alpha_{12}+\alpha_{12}\alpha_{13}-\frac{\alpha_{11}\alpha_{12}\alpha_{34}}{\alpha_{33}}}{\alpha_{44}}F^t_4;} \\

\\

F^t_1F^t_2 &=& \multicolumn{4}{l}{\frac{\alpha_{11}\alpha_{22}}{\alpha_{33}}F^t_3+\frac{\alpha_{11}\alpha_{23}+\alpha_{13}\alpha_{22}-\frac{\alpha_{11}\alpha_{22}\alpha_{34}}{\alpha_{33}}}{\alpha_{44}}F^t_4;} \\

\\

F^t_2F^t_1 &=& \frac{\alpha_{11}\alpha_{22}+\alpha_{12}\alpha_{23}}{\alpha_{44}}F^t_4; &
F^t_2F^t_2 &=& \frac{\alpha_{22}\alpha_{23}}{\alpha_{44}}F^t_4; \\
F^t_1F^t_3 &=& \frac{\alpha_{11}\alpha_{33}}{\alpha_{44}}F^t_4; &
F^t_3F^t_1 &=& \frac{\alpha_{12}\alpha_{33}}{\alpha_{44}}F^t_4; \\
F^t_3F^t_2 &=& \frac{\alpha_{22}\alpha_{33}}{\alpha_{44}}F^t_4.\\

\end{array}$$

To make the computations easier, we will consider a new basis $f_1,f_2,f_3,f_4$ in $\mathfrak{N}_{03}(\alpha)$ such that 
\[f_2f_3=0, f_3f_3=0 \mbox{  and }f_4\mathfrak{N}_{03}(\alpha)=\mathfrak{N}_{03}(\alpha)f_4=0.\] 
Such a basis can be defined as 
\[f_1=e_1, f_2=e_2, f_3=e_1+\alpha e_2+i\sqrt{\alpha^2+1}e_3 \mbox{ and }f_4=e_4.\] 
The multiplication table of $\mathfrak{N}_{03}(\alpha)$ with this new basis is

$$\begin{array}{rclrcl}

f_1f_1 &=& f_4; &
f_2f_2 &=& f_4; \\
f_1f_2 &=& \alpha f_4; &
f_2f_1 &=& -\alpha f_4; \\
f_1f_3 &=& (1+\alpha^2)f_4; &
f_3f_1 &=& (1-\alpha^2)f_4; \\
f_3f_2 &=& 2\alpha f_4.\\

\end{array}$$

Some routine calculations show that taking

$$\begin{array}{rclrclrclrcl}

\alpha_{11} &=& (1+\alpha^2)t; &
\alpha_{12} &=& (1-\alpha^2)t; &
\alpha_{13} &=& -2t; &
\alpha_{14} &=& 0; \\
 &  &   & 
\alpha_{22} &=& 2\alpha t; &
\alpha_{23} &=& -2\alpha t; &
\alpha_{24} &=& 0; \\
 &  &   & 
  &  &   & 
\alpha_{33} &=& -4\alpha^2 t; &
\alpha_{34} &=& -4\alpha^2 t\frac{\alpha^2 -3}{\alpha^2+1}; \\
 &  &   & 
  &  &   &
    &  &   &
\alpha_{44} &=& -4\alpha^2 t^2, \\
\end{array}$$

we obtain exactly

$$\begin{array}{rclrcl}
F^0_1F^0_1 &=& F^0_4; &
F^0_2F^0_2 &=& F^0_4; \\
F^0_1F^0_2 &=& \alpha F^0_4; &
F^0_2F^0_1 &=& -\alpha F^0_4; \\
F^0_1F^0_3 &=& (1+\alpha^2)F^0_4; &
F^0_3F^0_1 &=& (1-\alpha^2)F^0_4; \\
F^0_3F^0_2 &=& 2\alpha F^0_4.
\end{array}$$

Then, it suffices to take \[E_1^t=F^t_1, E_2^t=F^t_2, E_3^t=\frac{i}{\sqrt{\alpha^2+1}}(F^t_1+F^t_2-F^t_3)\mbox{ and }E_4^t=F^t_4\] so that we have the desired degeneration ${\mathcal B}^{4}_{10} \to {\mathfrak N}_{03}(\alpha)$, by the method described in the previous subsection. Namely,

\begin{align*}
E_1^t&=t((1+\alpha^2)e_1+(1-\alpha^2)e_2-2e_3), \\
E_2^t&=2\alpha t(e_2-e_3), \\
E_3^t&=\frac{it}{\sqrt{1+\alpha^2}}\Big((1+\alpha^2)e_1 + (1+\alpha^2)e_2 -  2(1-\alpha^2)e_3 + \frac{4\alpha^2(\alpha^2-3)}{1+\alpha^2}e_4\Big), \\
E_4^t&=-4\alpha^2 t^2e_4. \\
\end{align*}

Regarding the pure nilpotent bicommutative algebras, similar computations show that $\mathcal{B}_{10}$, $\mathcal{B}_{12}$, $\mathcal{B}_{14}$, $\mathcal{B}_{20}$ or ${\mathcal B}^4_{24}(\alpha)$ degenerate to them. The explicit degenerations can be seen in Table B (Appendix A).
\end{Proof}

         \section*{Appendix A.}

\[\begin{array}{|l|c|llllll|}

\multicolumn{8}{c}{ \mbox{ {\bf Table A.}
{\it The list of $4$-dimensional nilpotent ``pure'' bicommutative algebras.}}} \\
\multicolumn{8}{c}{ } \\
\hline

{\mathcal B} &  \mathfrak{Der}  & \multicolumn{6}{c|}{\mbox{Multiplication table} }\\ \hline
{\mathcal B}^4_{01}  & 6   & e_1 e_1 = e_2  & e_2 e_1=e_3 &&&& \\ \hline
{\mathcal B}^4_{02}(\alpha) & 6   & e_1 e_1 = e_2 & e_1 e_2=e_3 & e_2 e_1=\alpha e_3 &&  &\\ \hline

{\mathcal B}^4_{03}  & 4  & e_1 e_1 = e_2 & e_1 e_2=e_4 & e_2 e_1=e_3 &&  &\\ \hline
{\mathcal B}^4_{04}(\alpha)  & 4   &
e_1 e_1 = e_2 &   e_1e_2=e_4 & e_2e_1=\alpha e_4 & e_3e_3=e_4 & &\\ \hline

{\mathcal B}^4_{05}  & 4  &
e_1 e_1 = e_2 &   e_1e_2=e_4 & e_1e_3= e_4 & e_2e_1= e_4 & e_3e_3=e_4 &\\ \hline

{\mathcal B}^4_{06}(\alpha\neq0)  & 5   &
e_1 e_1 = e_2 &   e_1e_2=e_4 & e_1e_3=e_4 & e_2e_1=\alpha e_4&  &\\ \hline

{\mathcal B}^4_{07}  & 4   &
e_1 e_1 = e_2 &    e_2e_1= e_4 & e_3e_3=e_4&& &\\ \hline

{\mathcal B}^4_{08}  & 5   &
e_1 e_1 = e_2 & e_1e_3=e_4 &    e_2e_1= e_4 && &\\ \hline

{\mathcal B}^4_{09}  & 5   &
e_1 e_1 = e_2 & e_1e_2=e_4 &    e_3e_1= e_4 && &\\ \hline

{\mathcal B}^{4}_{10}  & 2   &     e_1 e_2=e_3  & e_1 e_3=e_4 & e_2e_1=e_4 & e_3e_2=e_4&&\\ \hline

{\mathcal B}^{4}_{11}  & 3   &   e_1e_2=e_3 & e_1 e_3=e_4  & e_3 e_2=e_4&& &\\ \hline

{\mathcal B}^{4}_{12}  & 3   &  e_1e_2=e_3 & e_1 e_1=e_4 & e_2 e_1=e_4 & e_3 e_2=e_4 &&\\ \hline

{\mathcal B}^{4}_{13}  & 4   &  e_1e_2=e_3 & e_2 e_1=e_4  & e_3 e_2=e_4 &&&\\ \hline

{\mathcal B}^{4}_{14}  & 3   &  e_1e_2=e_3 & e_1 e_3=e_4 & e_2 e_1=e_4  & e_2 e_2=e_4& &\\ \hline

{\mathcal B}^{4}_{15}  & 4   &  e_1e_2=e_3 & e_1 e_3=e_4 & e_2 e_1=e_4&& &\\ \hline

{\mathcal B}^{4}_{16}  & 4   &  e_1e_2=e_3 &  e_1 e_3=e_4 & e_2 e_2=e_4 && &\\ \hline

{\mathcal B}^{4}_{17}  & 5   &  e_1e_2=e_3 &  e_1 e_3=e_4  &&&&\\ \hline

{\mathcal B}^{4}_{18}  & 4   &  e_1e_2=e_3 &  e_1 e_1=e_4  & e_3 e_2=e_4 && &\\ \hline

{\mathcal B}^{4}_{19}  & 5   &  e_1e_2=e_3  & e_3 e_2=e_4 &&& &\\ \hline

{\mathcal B}^{4}_{20}  & 3   &  e_1e_1=e_2 & e_2 e_1= e_3 &  e_1 e_2=e_4 & e_3 e_1=e_4&  &\\ \hline

{\mathcal B}^4_{21}  & 4   &
e_1 e_1 = e_2 &  e_2 e_1=e_3 &  e_3 e_1=e_4&& &\\ \hline

{\mathcal B}^4_{22} & 3 &  e_1 e_1 = e_2,  & e_1 e_2=e_3,  & e_1e_3=e_4, & e_2 e_1 = e_4  &&\\
\hline

{\mathcal B}^4_{23} & 3 &  e_1 e_1 = e_2,  & e_1 e_2=e_3,  & e_1e_3=e_4, & e_2 e_1 =  e_3 +e_4, & e_2e_2= e_4, & e_3e_1= e_4  \\
\hline

{\mathcal B}^4_{24}(\alpha) & 3 \ (\alpha\neq 0,1) &  e_1 e_1 = e_2,  & e_1 e_2=e_3,  & e_1e_3=e_4, & e_2 e_1 = \alpha e_3, & e_2e_2=\alpha e_4, & e_3e_1=\alpha e_4  \\
\hline

\end{array}\]

 \

\begin{landscape} 
{\tiny 
$$
\begin{array}{|rcl|llll|}

\multicolumn{7}{c}{ \mbox{ {\bf Table B.}
\mbox{\it  Degenerations of bicommutative algebras of dimension $4$.}}} \\
\multicolumn{7}{c}{}\\

  \hline
  {\mathcal B}^{4}_{10} &\to& {\mathfrak N}^4_{02}(\alpha)  & E_1^t=t(e_1+e_3),&  E_2^t=-t e_1+t(1-\alpha)e_2,&  E_3^t=t^2e_4,&  E_4^t= t^2(1-\alpha) (e_3+e_4) \\
  \hline
  {\mathcal B}^{4}_{10} &\to& {\mathcal B}^{4}_{01}  & E_1^t=te_1+e_2,  &E_2^t=t(e_3+e_4),  &E_3^t=te_4,  &E_4^t=te_2 \\
  \hline
  {\mathcal B}^{4}_{10} &\to& {\mathcal B}^{4}_{02}(\alpha)_{\alpha\neq0}  & E_1^t=e_1+\alpha e_2,&  E_2^t=\alpha(e_3+e_4),&  E_3^t=\alpha e_4,&  E_4^t=t(e_2+e_3) \\
  \hline
  {\mathcal B}^{4}_{20} &\to& {\mathcal B}^{4}_{03}  & E_1^t=te_1,&  E_2^t=t^2e_2,&  E_3^t=t^3e_3,&  E_4^t=t^3e_4 \\
  \hline
  {\mathcal B}^{4}_{10} &\to& {\mathcal B}^4_{04}(\alpha)  & E_1^t=-t^2e_1-\alpha t^2e_2+((\alpha+1)t^2+t^4)e_3,&  E_2^t=\alpha t^4e_3+(t^4(\alpha-(\alpha+1)(\alpha+1+t^2)))e_4,&  E_3^t=t^3e_1-\alpha t^3e_3,&  E_4^t=-\alpha t^6 e_4 \\
  \hline
  {\mathcal B}^{4}_{10} &\to& {\mathcal B}^{4}_{05}  & E_1^t=t^2(e_1+e_2+(it-2)e_3),&  E_2^t=t^4(e_3+(2it-3)e_4),&  E_3^t=it^3(e_2-e_3),&  E_4^t=t^6e_4 \\
  \hline
  {\mathcal B}^{4}_{10} &\to& {\mathcal B}^4_{06}(\alpha)_{\alpha\neq0}  & E_1^t=e_1+\alpha e_2-\alpha(1+\alpha t^{-1})e_3,&  E_2^t=\alpha e_3-\alpha^2(1+(1+\alpha)t^{-1})e_4,&  E_3^t=te_2,&  E_4^t=\alpha e_4 \\
  \hline
  {\mathcal B}^{4}_{10} &\to& {\mathcal B}^{4}_{07}  & E_1^t=t^4e_1-t^2e_2+(t^2-t^4+t^6)e_3,&  E_2^t=-t^6e_3+(-t^4+t^6-2t^8+t^{10})e_4,&  E_3^t=t^{5}e_1+t^{3}e_3,&  E_4^t=t^8e_4 \\
  \hline
  {\mathcal B}^{4}_{10} &\to& {\mathcal B}^{4}_{08}  & E_1^t=te_1+e_2-(1+t^{-1})e_3,&  E_2^t=te_3-(2+t^{-1})e_4,&  E_3^t=te_2+t^2e_3,&  E_4^t=te_4 \\
  \hline
  {\mathcal B}^{4}_{20} &\to& {\mathcal B}^{4}_{09}  & E_1^t=te_1,&  E_2^t=t^2e_2,&  E_3^t=t^2e_3,&  E_4^t=t^3e_4 \\
  \hline
  {\mathcal B}^{4}_{10} &\to& {\mathcal B}^{4}_{11}  & E_1^t=t^{-1}e_1,&  E_2^t=t^{-1}e_2,&  E_3^t=t^{-2}e_3,&  E_4^t=t^{-3}e_4 \\
  \hline
  {\mathcal B}^{4}_{10} &\to& {\mathcal B}^{4}_{12}  & E_1^t=te_1+e_3,  &E_2^t=e_2,  &E_3^t=te_3+e_4,  &E_4^t=te_4 \\
  \hline
  {\mathcal B}^{4}_{10} &\to& {\mathcal B}^{4}_{13}  & E_1^t=te_1,&  E_2^t=e_2,&  E_3^t=te_3,&  E_4^t=te_4 \\
  \hline
  {\mathcal B}^{4}_{10} &\to& {\mathcal B}^{4}_{14}  & E_1^t=e_1,  &E_2^t=te_2+e_3,  &E_3^t=te_3+e_4,  &E_4^t=te_4 \\
  \hline
  {\mathcal B}^{4}_{10} &\to& {\mathcal B}^{4}_{15}  & E_1^t=e_1,&  E_2^t=te_2,&  E_3^t=te_3,&  E_4^t=te_4 \\
  \hline
  {\mathcal B}^{4}_{14} &\to& {\mathcal B}^{4}_{16}  & E_1^t=t^{-1}e_1,&  E_2^t=t^{-2}e_2,&  E_3^t=t^{-3}e_3,&  E_4^t=t^{-4}e_4 \\
  \hline
  {\mathcal B}^{4}_{10} &\to& {\mathcal B}^{4}_{17}  & E_1^t=t^{-1}e_1,&  E_2^t=e_2,&  E_3^t=t^{-1}e_3,&  E_4^t=t^{-2}e_4 \\
  \hline
  {\mathcal B}^{4}_{12} &\to& {\mathcal B}^{4}_{18}  & E_1^t=t^{-2}e_1,&  E_2^t=t^{-1}e_2,&  E_3^t=t^{-3}e_3,&  E_4^t=t^{-4}e_4 \\
  \hline
  {\mathcal B}^{4}_{10} &\to& {\mathcal B}^{4}_{19}  & E_1^t=e_1,&  E_2^t=t^{-1}e_2,&  E_3^t=t^{-1}e_3,&  E_4^t=t^{-2}e_4 \\
  \hline
  {\mathcal B}^{4}_{24}(\frac{1}{t}) &\to& {\mathcal B}^{4}_{20}  & E_1^t=te_1+\frac{t}{1-t}e_2,&  
  E_2^t=t^2e_2+(1+t)\frac{t}{1-t}e_3+\frac{t}{(1-t)^2} e_4,&  E_3^t=t^2e_3+(1+2t)\frac{t}{1-t}e_4,&  E_4^t=t^2e_4 \\
  \hline
  {\mathcal B}^{4}_{20} &\to& {\mathcal B}^{4}_{21}  & E_1^t=t^{-1}e_1,&  E_2^t=t^{-2}e_2,&  E_3^t=t^{-3}e_3,& E_4^t=t^{-4}e_4 \\
  \hline
  {\mathcal B}^{4}_{24}(t) &\to& {\mathcal B}^{4}_{22}  & E_1^t=te_1+te_2,&  
  E_2^t=t^2e_2+(t^2+t^3)e_3+t^3 e_4,&  E_3^t=t^3e_3+(t^3+2t^4)e_4,&  E_4^t=t^4e_4 \\
  \hline
  {\mathcal B}^{4}_{24}(1-t) &\to& {\mathcal B}^{4}_{23}  & E_1^t=te_1+te_2,&  
  E_2^t=t^2e_2+(2t^2-t^3)e_3+t^2(1-t)e_4,&  E_3^t=t^3e_3+(3t^3-2t^4)e_4,&  E_4^t=t^4e_4 \\
  \hline

\end{array}$$
}

\end{landscape}


\end{document}